\RequirePackage{fix-cm}
\documentclass[smallextended]{svjour3}       
\smartqed  
\usepackage{graphicx}

\usepackage{amsmath}
\usepackage{amssymb}
\usepackage{epsfig}
\usepackage{colordvi}
\usepackage{graphics}
\usepackage{color}

\usepackage{cite}
\usepackage{rotating}
\usepackage{booktabs}
\usepackage{lmodern}

 \usepackage{mathptmx}      
 \journalname{Journal of Elasticity}

\newcommand{\bfb}{{\bf b}}%
\newcommand{\bfn}{{\bf n}}%
\newcommand{\bft}{{\bf t}}%
\newcommand{\bfvartheta}{\boldsymbol{\vartheta}}%
\newfont{\tenbfit}{cmbx10}%

\newfont{\tenbbb}{msbm10}%
\newfont{\svnbbb}{msbm8}%

\newcommand{\lj}{\mbox{$[\kern-0.1478125em[$}}
\newcommand{\rj}{\mbox{$]\kern-0.1478125em]$}}
\newcommand{\la}{\mbox{$\langle\kern-0.2325em\langle$}}
\newcommand{\ra}{\mbox{$\rangle\kern-0.2325em\rangle$}}
\newcommand{\Blj}{\mbox{$\Big[\kern-0.275em\Big[$}}
\newcommand{\Brj}{\mbox{$\Big]\kern-0.275em\Big]$}}
\newcommand{\Bla}{\mbox{$\Big\langle\kern-0.425em\Big\langle$}}
\newcommand{\Bra}{\mbox{$\Big\rangle\kern-0.425em\Big\rangle$}}

\newcommand{\bfmfrakK}{{\boldsymbol{\mathfrak{K}}}}
\newcommand{\bfmfrakM}{{\boldsymbol{\mathfrak{M}}}}

\begin{document}

\title{Translation of Michael Sadowsky's paper ``The differential equations of the {\sc{M\"obius}} band"\footnote{Citations of this translation should refer also to Sadowsky's original paper, as cited in the Abstract.}
}


\author{Denis F.\ Hinz         \and
        Eliot Fried 
}


\institute{Denis F. Hinz \at
              Mathematical Soft Matter Unit\\
              Okinawa Institute of Science and Technology \\
              Okinawa, Japan 904-0495\\
              \email{dfhinz@gmail.com}   
           \and
           Eliot Fried \at
              Mathematical Soft Matter Unit\\
              Okinawa Institute of Science and Technology, \\
              Okinawa, Japan 904-0495\\
              \email{eliot.fried@oist.jp}
}

\date{Received: date / Accepted: date}

\maketitle

\begin{abstract}
This article is a translation of Michael Sadowsky's original paper ``Die Differentialgleichungen des {\sc{M\"obius}}schen Bandes." in Jahresbericht der Deutschen Mathematiker-Vereini\-gung {\bf{39}} (2.\ Abt.\ Heft 5/8, Jahresversammlung vom 16.\ bis 23.\ September), 49--51  (1929), which is a short version of his paper ``Theorie der elastisch biegsamen undehnbaren B\"ander mit Anwendungen auf das {\sc{M\"obius}}'sche Band" in 3.\ internationaler Kongress f\"ur technische Mechanik,  Stockholm, 1930. 
\end{abstract}


\section*{Translation of the original paper}

In a previously completed work, Sadowsky~\cite{Sadowsky1930} provided a proof for the {\emph{existence}} of a developable {\sc{M\"obius}} band and provided a variational context for the underlying geometric problem. In the present work, the problem of determining the equilibrium shape of a {\sc{M\"obius}} band formed from an elastic strip is treated as a static mechanical problem.

Let $P$ denote an arbitrary point on the midline of the {\sc{M\"obius}} band, let $s$ denote the arclength of the mid-line, let $\bft$, $\bfn$, and $\bfb$ denote the accompanying vector triad at $P$, and let $K$ and $W$ denote the curvature and twist of the midline at $P$. Considering a cut perpendicular to the midline at $P$, the stresses along the line of the cut may be reduced to a force $\bfmfrakK$ and a moment $\bfmfrakM$ both acting at $P$.

The representations of $\bfmfrakK$ and $\bfmfrakM$ relative to the triad $\bft$, $\bfn$, and $\bfb$ shall be denoted as
\begin{equation}\label{eq:K}
\bfmfrakK = T\bft + N\bfn+ B\bfb,
\end{equation}
and
\begin{equation}\label{eq:M}
\bfmfrakM = {\mathfrak{T}} \bft + {\mathfrak{N}} \bfn + {\mathfrak{B}} \bfb.
\end{equation}
From a mechanical perspective, the band is defined through the following three requirements:
\begin{enumerate}
\item The band is an object with distinguished midline.
\item The vector triad $\bft$, $\bfn$, and $\bfb$ accompanying the midline is material.
\item The midline is inextensible.
\end{enumerate}
(The second requirement results from the band being part of the rectifying surface of its midline).

Applying a virtual deformation to the band and invoking the Frenet--Serret formulas results in
\begin{equation}\label{eq:virt01}
d \delta \bfvartheta = (\bft \delta W+ \bfb \delta K)ds,
\end{equation}
and in
\begin{equation}\label{eq:virt02}
 \delta d s  = 0,
\end{equation}
where $\delta \bfvartheta$ is the virtual rotation of the vector triad at $P$.

{\emph{The following considerations hold only for infinitesimally narrow bands.}}

The virtual work of the internal forces is\footnote{\textsc{Hamel}, \"Uber die Mechanik der Dr\"ahte und Seile, Sitzungsberichte der Berliner Mathematischen Gesellschaft, XXV, 1925/26. Also \textsc{Hamel} Die Axiome der Mechanik, Handbuch der Physik, herausgegeben von \textsc{H.\ Geiger} und \textsc{Karl Scheel}, Band V. In the latter work further references may be found.}
\begin{equation}\label{eq:dA}
\delta A_i = - \int \limits_{s_1}^{s_2} ({\mathfrak{T}} \delta W + {\mathfrak{B}} \delta K  )ds.
\end{equation}

The mean curvature $H$ of the band surface at $P$ is$^1$
\begin{equation}
H = \frac{K^2+W^2}{2K}.
\end{equation}

The elastic potential $E$ of the band is proportional to the square of the mean curvature (granted the presumption that the band surface is a developable surface)
\begin{equation}\label{eq:elast}
E = A \frac{(K^2+W^2)^2}{K^2}.
\end{equation}

Since $\delta A_i = - \int \limits_{s_1}^{s_2} \delta E d s$ and, in view of~\eqref{eq:dA} and~\eqref{eq:elast}, it follows that
\begin{equation}\label{eq:08}
{\mathfrak{T}} = \frac{\partial E}{\partial W} = A \frac{4W(K^2+W^2)}{K^2}
\end{equation}
and
\begin{equation}\label{eq:09}
{\mathfrak{B}} = \frac{\partial E}{\partial K} = A \frac{2(K^4-W^4)}{K^3},
\end{equation}
where $A$ is a positive material constant. The six components appearing in~\eqref{eq:K} and~\eqref{eq:M} satisfy six known {\emph{equilibrium equations}}.$^1$ Using \eqref{eq:08} and~\eqref{eq:09}, those equations read

\begin{equation}\label{eq:relation02}
\left\{
\begin{array}{l}
\displaystyle
T =  A\left( C - \frac{(K^2+W^2)^2}{K^2} \right),
\cr\noalign{\vskip8pt}
\displaystyle
N = - \frac{A}{K} \frac{d}{ds} \frac{(K^2+W^2)^2}{K^2},
\cr\noalign{\vskip8pt}
\displaystyle
B = \frac{2AW}{K^3}(K^2+W^2)^2 +4A \frac{d}{ds} \left( \frac{1}{K} \frac{d}{ds} \frac{W (K^2+W^2)}{K^2}\right),
\cr\noalign{\vskip8pt}
\displaystyle
\mathfrak N = \frac{4A}{K}\frac{d}{ds} \frac{W (K^2+W^2)}{K^2},
\end{array}
\!\!
\right.
\end{equation}
and 
\begin{equation}\label{eq:relation03}
\left\{
\begin{array}{l}
\displaystyle
KT + \frac{dN}{ds} - WB = 0,
\cr\noalign{\vskip8pt}
\displaystyle
WN + \frac{dB}{ds} = 0,
\end{array}
\!\!
\right.
\end{equation}
where $C$ is an arbitrary constant of integration. 

If the expressions~\eqref{eq:relation02} were to be used in~\eqref{eq:relation03}, two differential equations for $K$ and $W$, respectively, would result. The integration of these equations would yield the natural midline of the band. Equations~\eqref{eq:08}, \eqref{eq:09}, \eqref{eq:relation02}, and~\eqref{eq:relation03} shall be referred to as the equations of the {\sc{M\"obius}} band.

The midline of the band possesses a singular point $X$ that can be found in the following way: a {\sc M\"obius} band possesses a symmetry axis defined such that it is congruent with itself after a rotation of $180^\circ$ about that axis. The symmetry axis intersects the band at two points such that the axis coincides with the binormal $\bfb$ at one of these points; this is the point $X$.

Let $\varphi$ denote the angle between the rectilinear generators of the band through the point $P$ and $\bfb$; then$^1$
\begin{equation}
\tan \varphi = \frac{W}{K}.
\end{equation}

For the point $X$,
\begin{equation}\label{eq:Beq0}
{\mathfrak B} =0
\end{equation}
due to symmetry.

To arrive at an additional conclusion, consider an {\emph{experiment with a band model}}: this shows that
\begin{equation}\label{eq:exp01}
\lim \limits_{P \rightarrow X} \varphi \neq 0,
\end{equation}
\begin{equation}\label{eq:exp02}
\lim \limits_{P \rightarrow X} {\mathfrak T} \neq 0.
\end{equation}
The experiment corresponding to~\eqref{eq:exp01} consists of observing a band; for~\eqref{eq:exp02}, one needs to cut the boundaries of a band at $X$ to infer the moment $\mathfrak T$ from the twist of the band in the weakened cross-section. From~\eqref{eq:Beq0},~\eqref{eq:exp01}, \eqref{eq:exp02}, and the governing equations~\eqref{eq:08} and~\eqref{eq:09} of the band it transpires that
\begin{equation}\label{eq:lim01}
\lim \limits_{P \rightarrow X} K \neq 0, \quad \lim \limits_{P \rightarrow X} W \neq 0, \quad \text{and} \quad \lim \limits_{P \rightarrow X} \varphi = 45^\circ.
\end{equation}
The last limit in~\eqref{eq:lim01} points to the following peculiar fact:\\

{\emph{A {\sc M\"obius} band consists of a planar, right triangle. The curved, analytic portion of the band connects continuously with the two legs of the right triangle; it connects with continuous tangential plane, but with discontinuous curvature.\footnote{For an interpretation of the last section of Sadowsky's original paper including figures generated from recent numerical simulations, refer to D.~F.~Hinz and E.~Fried, Translation and interpretation of Michael Sadowsky's paper ``Theory of elastically bendable inextensible bands with applications to the {\sc{M\"obius}} band", submitted to \emph{Journal of Elasticity} (2014). }}


\begin{thebibliography}{9}
%
%



\bibitem{Sadowsky1930}
M.~Sadowsky, Ein elementarer Beweis f\"ur die Existenz eines abwickelbaren M\"obiusschen Bandes und die Zur\"uckf\"uhrung des geometrischen Problems auf ein Variationsproblem, Sitzungsberichte der Preussischen Akademie der Wissenschaften, physikalisch-mathematische Klasse {\bf 22}, 412--415 (1930)

%
%
%




\end{thebibliography}
\end{document}